\documentclass{ifacconf}

\usepackage{graphicx}      % include this line if your document contains figures
\usepackage{natbib}        % required for bibliography
\let\theoremstyle\relax
\usepackage{amsthm}
\usepackage{amsfonts}
\usepackage{amsmath}
\usepackage{graphicx} 
\usepackage{epsfig}
\usepackage{bbm}
\usepackage{epstopdf}
\newtheorem{proposition}{Proposition}

\newtheorem{remark}{Remark}

\DeclareMathOperator{\PP}{\mathbb{P}}
\DeclareMathOperator{\Q}{\mathbb{Q}}
\DeclareMathOperator{\A}{\mathcal{A}}
\DeclareMathOperator{\T}{\mathcal{T}}

\DeclareMathOperator{\R}{\mathbb{R}}

\newcommand{\I}{{\mathcal I}}

\newcommand{\RR}{{\mathbb R}}

\newcommand{\lra}{\longrightarrow}

\newcommand{\cI}{{\mathcal I}}

\newcommand{\cA}{{\mathcal A}}

\newcommand{\cN}{{\mathcal N}}

\newcommand{\cP}{{\mathcal P}}

\newcommand{\cT}{{\mathcal T}}

\newcommand{\cQ}{{\mathcal Q}}

\renewcommand{\phi}{\varphi}

\begin{document}
\begin{frontmatter}

\title{The orienteering problem: a hybrid control formulation} 

\thanks[footnoteinfo]{The second author was supported by MIUR grant ``Dipartimenti Eccellenza 2018-2022" CUP: E11G18000350001, DISMA, Politecnico di Torino}

\author[Bagagiolo]{Fabio Bagagiolo} 
\author[Festa]{Adriano Festa}
\author[Marzufero]{Luciano Marzufero}

\address[Bagagiolo]{Dipartimento di Matematica,
Universit\`a di Trento, 
Via Sommarive, 14, 38123 Povo (TN) Italy, (e-mail: fabio.bagagiolo@unitn.it).}
\address[Festa]{Dipartimento di Scienze Matematiche ``G. L. Lagrange", 
Politecnico di Torino, 
Corso Duca degli Abruzzi, 24, 10129 Torino Italy, (e-mail: adriano.festa@polito.it)}
\address[Marzufero]{Dipartimento di Matematica,
Universit\`a di Trento, 
Via Sommarive, 14, 38123 Povo (TN) Italy, (e-mail: luciano.marzufero@unitn.it).}

\begin{abstract}                % Abstract of not more than 250 words.
In the last years, a growing number of challenging applications in navigation, logistics, and tourism were modeled as orienteering problems. This problem has been proposed in relation to a sport race where certain control points must be visited in a minimal time. In a certain kind of these competitions, the choice of the number and the order for the control points are left to the competitor. We propose an original approach to solve the orienteering problem based on hybrid control. The continuous state of the system contains information about the navigation toward the next control point. In contrast, the discrete state keeps track of the already visited ones. The control problem is solved using non-standard dynamical programming techniques.
\end{abstract}


\end{frontmatter}
%===============================================================================

\section{Introduction}

\begin{figure}[t]
\begin{center}
\begin{tabular}{c}
\includegraphics[width=8.6cm]{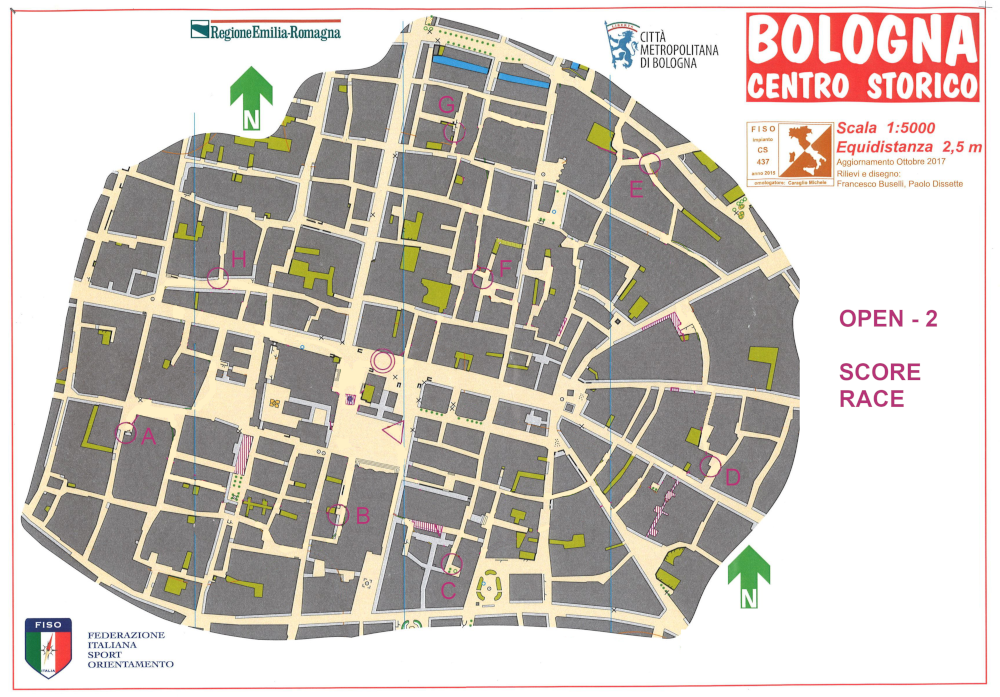}
\end{tabular}
\end{center}
\caption{An example of an urban race taking place in the historical city center of Bologna, IT. In the notation the start is marked by a triangle, the checkmarks by circles and the arrival by the double circle. \emph{(credits to Open Maps Project of OrienteeringTime.com)}}\label{1}
\end{figure}

The ``Orienteering Problem'' originates from the sport of orienteering \cite{chao1996fast}. This sports game is performed in various kinds of environments as parks, rural or mountain areas, but also inside the cities.  The competitors, with the help of a detailed map (see Fig. \ref{1}) and a compass, which are the only technological instruments allowed during the race, start at a specific control point (\emph{start}, normally marked with a triangle), must visit a collection of checkpoints (circles) as fast as it is possible, to return to a final control point (\emph{arrival} normally marked with a double circle). Depending on the competition type, the checkpoints could be visited in a specific order (\emph{cross-country}) or an order decided by the competitor (\emph{free-order}). The placings are determined, in general, by the total time of the race of each player.  In the special \emph{score} competition type, at each of the control points is associated a score, and the competitor must collect as many points as possible in a given time, including the return to the arrival. 

As a consequence, the orienteering problem includes elements of optimal navigation through inhomogeneous media as minimal time control problems, within the classic Traveling Salesmen problems \cite{bellman1962dynamic, gavish1978travelling} and Knapsack problems \cite{martello1987algorithms}. The mathematical literature about the problem includes mainly combinatorial optimization approaches (for a survey about these contributions see \cite{vansteenwegen2011orienteering}, or the more recent \cite{gunawan2016orienteering}). Generally, the problem is modeled on a graph, not considering its navigation component. 

This paper wants to propose an original approach to the complete Orienteering Problem based on a hybrid control formulation using dynamical programming techniques and semiLagrangian numerical approximation. In particular, we model the player's strategic choice as the optimization of a trajectory across a given 2D domain where an additional discrete state records the visited checkpoints. This model shares various similarities with the ``Optimal Visiting Problem'', as proposed in \cite{bagfagmagpes}, and also discussed in \cite{BFMProc1, BFMNoDea}, where other aspects of the framework are explored.  

The paper is organized as follows. Section \ref{s:hyb} defines the model, in its stationary formulation, that will be used in the cross-country and free cases. The time-dependent variant, useful in the score case, is sketched in the same section. Section \ref{s:num} discusses the numerical techniques used to approximate the Hamilton-Jacobi-Bellman (HJB) equation, which describes the value function of the optimal control problem. In Section \ref{s:test}, we show various tests that display the excellent qualities of our technique.  Section \ref{s:con} concludes the paper and points out some interesting open research questions.

\section{A hybrid control formulation of the orienteering problem}\label{s:hyb}

We propose to model the orienteering problem as an hybrid control problem. Given a collection of  $m$  distinct checkpoints $\{\T_j\}_{j=1,\ldots,m}\subset \R^2$ and an arrival $\A\in\R^2$, we identify the finite number of states of the system as $\I=\{1,2,\ldots,\hat q\}$ isomorphic to a subset of all possible permutations with repetitions of the elements $\{0,1\}$ where $0$ means that the $i$-checkpoint has not been visited and viceversa. Therefore, at the start the player is in the state $1\cong (0,\ldots,0)=(\hat q_1,\ldots,\hat q_{m})$, while in $\hat q\cong (1,\ldots,1)$ the competitor has already visited all the control points and has only to reach the arrival. 

The controlled dynamics are described by:
\begin{equation} \label{eq_stato}
 \begin{cases}
 y'(t)=f(y(t),\alpha(t),Q(t))\\
 y(0)= x, \ Q(0^+)=q,
 \end{cases}
\end{equation}
where $x,y\in \R^2$, $Q,q\in \I$. The function $f :\R^2\times A \times \I\lra \R^2$ represents the continuous dynamics, globally bounded and uniformly Lipschitz continuous w.r.t.the space variable and $\alpha:[0,\infty) \rightarrow A$ is a control function.

The term $Q(t)$ models the possibility to switch between the various dynamics of the system with different targets and takes values in the set of piecewise constant discrete controls $\cQ$, that is:
\begin{equation}\label{switch}
\cQ = \left\{Q(\cdot):[0,\infty) \lra \I \> | \> Q(t)=\sum_{i=0}^m w_i \chi_{t_i}(t)\right\},
\end{equation}
where $\chi_i(t)=1$ if $t\in [t_i,t_{i+1})$ and $0$ otherwise,  $\{t_i\}_{i=0,\ldots,m}$ are the (ordered) times at which a switch occurs (with the convention that $t_0=0$), and $\{w_i\}_{i=1,\ldots,m}$ are values in $\I$. The maximum number of switches is limited by $m$ (the number of checkpoints) even if the number of states is in general bigger. This is due, as we describe below, by the design of the switching cost. Intuitively, a player cannot discard more than the total number of the checkpoints. Clearly, in the case of a multiple switch (simultaneous renounce to more than one checkpoint), some of these terms may coincide. 
The choice of the control strategy defined as $\mathcal{S}:=\left(\alpha,Q\right)$ has the objective of minimizing the following cost functional of minimum time type:
\begin{multline}\label{J}
 J(x,q;\mathcal{S})  := \int_{t_{m}}^{\tau(x,q)}\ell(y_{x,q}(s))e^{-\lambda s} ds  \\+ \sum_{i=1}^{m-1}\int_{t_i}^{t^{i+1}}\ell(y_{x,q}(s))e^{-\lambda s} ds \\
 + \sum_{i=0}^{m} e^{-\lambda t_i}C\left(y_{x,q}(t_i),Q(t_i^-),Q(t_i^+)\right),
\end{multline}
where  $\lambda>0$ is a discount factor, $\ell$ is a running cost function and the trajectory $y_{\alpha,t,x}$ is a solution of \eqref{eq_stato} starting at time $t=0$ in the point $(x,q)$. The term $\tau(x,q)$ is the first time of arrival on the arrival $\A$ starting from $(x,q)$.

In order to have well position of the model, the switching cost $C:\R^2\times \I \times \I  \lra \R_+$ should have a strictly positive infimum, to be bounded and Lipschitz continuous w.r.t.~$x$ (see for details \cite{BFMProc1}). Since we want that the switches happen only very close to the checkpoints we model the switch cost in a penalization fashion as 
\begin{equation}
C(x,q,q')=\
\frac{1}{\varepsilon}\sum_{j\in  \cP(q)}\|x-\cT_j\|, \\
\end{equation}
where the map $\cP:\cI \longrightarrow P \subseteq \{0,1,\ldots,2^m\}$ is a multivalued function which indicates the permitted switching and it can be chosen as:
\begin{enumerate}
\item[a)]\emph{Cross-country} $\cP(q)=\{q+1\}$ and  $\cP(\hat q)=\emptyset$.
\item[b)]\emph{Free-order} $$\cP(q)=\{q'\in \I\, |\, q'_i=1 \hbox{ for every }i \hbox{ s.t. }q_i=1\}.$$
\end{enumerate}

The \emph{value function} $v$ of the problem is then defined, for $\mathcal{S}\in\cA\times\cQ$, as:
\begin{eqnarray} \label{f_valore}
 v(x,q):= \inf_\mathcal{S} J(x,q;\mathcal{S}),
\end{eqnarray}
and is characterized via a suitable Hamilton-Jacobi-Bellman (HJB) equation. Continuity of the value function, which is a typical result in the more classic theory of viscosity solutions is not an easy task in deterministic hybrid control problems (for a precise statement of the hypothesis we refer to \cite{BFMProc1} or to \cite{dharmatti2005hybrid}), even though the most general results have been proved in a weaker framework (see \cite{bensoussan1997hybrid}). 

In the present framework the \emph{Dynamical Programming Principle} (DPP) holds (for the proof we refer to \cite{dolcetta1984optimal}).

\begin{proposition}
For each $\{t_i\}_{i=1,\ldots,2^m}$ there exists a $\hat q\in\I$ such that
\begin{equation}\label{strict}
   v(x,\hat q)<\min\limits_{q\neq\hat q}(v(x,q)+C(x,\hat q,q)).
\end{equation}
Moreover, for any $\hat t\in [t_i, t_{i+1})$, we have 
$$v(x,\hat q)=\int_{\hat t}^{t}\ell(y_{x,q}(s))e^{-\lambda s} ds+e^{-\lambda t}v(y_x(t),\hat q).
$$
Finally, as a consequence, the following DPP holds
\begin{multline}\label{DPP}
 v(x,\hat q)=\int_0^{t}\ell(y_{x,q}(s))e^{-\lambda s} ds 	\\
 + \sum_{i=0}^K e^{-\lambda t_i}C\left(y_{x,q}(t_i),Q(t_i^-),Q(t_i^+)\right)\\+v(y_x(t),Q(t_K^+)),
\end{multline}
for any $t\in[0,\tau(x,q)]$, and $K=\max_i\{t_i<t\}$.
\end{proposition}

Using the DPP above for a small time step, we can prove that the value function of the problem solves a HJB equation in a Quasi-Variational Inequality form. In other words, defining for $x,p\in\R$ and $i\in\I$ the Hamiltonian function by
\begin{equation}\label{Ham}
 H(x,q,p) := \sup_{\alpha\in A}\{ - f(x,\alpha,q)\cdot p -\ell(x)) \}
\end{equation}
and the switching operator $\cN$ by:
\begin{equation*}
 \cN \phi(x,\alpha) := \inf_{\beta\in \I} \{\phi(x,\beta)+C(x,\alpha,\beta)\},
\end{equation*}
for every function $\phi:\R^2\times \I\lra\R_+$, we have the following result (for a proof of a similar case refer to \cite{dharmatti2005hybrid, ferretti2019optimal}).

\begin{proposition}
Under the hypothesis of the current section the value function defined in \eqref{f_valore} is the unique viscosity solution of
  \begin{multline}\label{hjb}
 \max\left(v(x,q)-\cN v(x,q), \right. \\ \left.\hspace{2cm}\lambda v(x,q)+H(x,q,D v(x,q)) \right)= 0,  
\end{multline}
for any $(x,q)\in (\R^2\times \cI)\setminus(\cA\times\{1\})$ and with $v(x,1)=0, $  if $x\in\cA$.
\end{proposition}

The equation \eqref{hjb} is HJB in the special form of a system of Quasi-Variational Inequalities complemented with some boundary condition in the typical form of a \emph{minimal time  optimal control} problem.

In \eqref{hjb} there are contained two separate Bellman operators which provide respectively the best possible switching, and the best possible continuous control. The argument attaining the maximum in \eqref{hjb} represents the overall optimal control strategy.

\begin{remark}
The cross-country case (i.e. the order of the switches between states is fixed) trivializes the problem, which degenerates to a collection of disjoint minimal time problems that may be solved separately. Our model, anyway, shares with the latter approach the same degree of complexity. 
\end{remark}

\subsection{A time-dependent case: the \emph{score} competition}
As stated before, in the score competition the player has to keep his race in a fixed time $T>0$ visiting as many checkpoints as it is possible. This changes the model that has to be time dependent. Actually, as shown in \cite{Bardi, Festa2017127}, the time dependent case can be seen as a special case of the stationary one, through a duplication of variables procedure. Therefore, the theoretical results valid for the stationary case, can be generally adapted to this new case. We limit ourself to report the form of the cost functional that is
\begin{multline}\label{Jt}
 J(t,x,q;\mathcal{S})  :=  \sum_{i=0}^{k}\int_{t_i}^{t^{i+1}}e^{-\lambda (s-t)} ds \\
 + \sum_{i=1}^{k} e^{-\lambda (t_i-t)}C\left(y_{x,q}(t_i),Q(t_i^-),Q(t_i^+)\right),
\end{multline}
where all the switches happen in the interval $[t,T]$, i.e. $t=t_0\leq t_1\leq\ldots\leq t_{k}\leq T$. The value function 
$$ v(t,x,q):=\inf_{S} J(t,x,q;\mathcal{S})$$
 is characterized as viscosity solution of the HJB equation in the variational inequality form as
  \begin{multline}\label{hjbt}
 \max\left(v(t,x,q)-\cN v(t,x,q), \right. \\ \left.\hspace{0.3cm}-v_t(t,x,q)+\lambda v(t,x,q)+H(x,q,D v(t,x,q)) \right)= 0,  
\end{multline}
for any $(t,x,q)\in (\R^2\times \cI)\setminus(\cA\times\{1\})$ and with $v(t,x,1)=0, $  if $x\in\cA$ and for any $t\leq T$.
A more detailed presentation of the theoretical results in the time-dependent case is \cite{BFMProc1}.

\section{Numerical approximation}\label{s:num}

In order to set up a numerical approximation for \eqref{hjb}, we consider a discrete grid of nodes $(x_{i},q)$ in the state space $\Omega$  with discretization parameters $\Delta x$ that for simplicity we consider as the open set $[a,b]\times[c,d]$, with $a,b,c,d \in\R$. Thus, a general $i=(i_1,i_2)$ $x_{i}=(a+i_1 \Delta x , c+i_2\Delta x)$. We introduce an auxiliary parameter $h$ which has the role of tuning the distance of the search of information in the upwind direction and it should be taken, generally, of order $O(\Delta x)$ to maximize the accuracy of the approximation \cite{FF13}. In what follows, we denote the discretization steps in compact form by $\Delta=(h, \Delta x)$ and the approximate value function by $V^\Delta$.

Following \cite{MR3328207}, we write the scheme at $(x_j,q)$ in fixed point form as
\begin{equation}\label{Scheme1}
V_\Delta(x_j,q) = \min \left( NV_\Delta(x_j,q), \Sigma\left(x_j,q,V_\Delta\right)\right).
\end{equation}
In \eqref{Scheme1}, the numerical operator $\Sigma$ is related to the continuous control, or, in other terms, to the approximation of the Hamiltonian function \eqref{Ham}. The discrete switch operator $N$ is computed at a node $(x_j,q)$ as
\begin{equation}\label{Scheme11}
NV_\Delta(x_j,q) := \min_{q'\in \I}\left\{V_\Delta(x_j,q') + C(x_j,q,q')\right\},
\end{equation}
and, in fact, this corresponds to the exact definition.

\subsection{A Semi-Lagrangian scheme.} A viable technique for solving \eqref{hjb} is a semi-Lagrangian scheme, obtained by adapting the scheme proposed in \cite{camilli1995approximation}. The main advantage of such approach is an unconditional stability of the scheme with respect to the discretization parameters, still keeping monotonicity.

The scheme requires extending the node values to all $x\in\RR^d$ using an interpolation $\mathbb{I}$. We denote by $\mathbb{I}\left[V_\Delta\right] (x,q)$ the interpolation of the values $V_\Delta(x_j,q)$ computed at $(x,q)$. With this notation, a standard semi-Lagrangian discretization of the Hamiltonian is given (see \cite{camilli1995approximation}) by
\begin{multline}\label{Scheme2}
\Sigma\left(x_j,q,V_\Delta\right)= \\
h\ell(x_j) + e^{-\lambda h}\min_{u\in U}\left\{\mathbb{I}\left[V_\Delta\right] (x_j+h \> f(x_j,q,u), q)\right\}.
\end{multline}
The full scheme is obtained by using \eqref{Scheme11}-\eqref{Scheme2} in \eqref{Scheme1}.

As far as the interpolation $\mathbb I$ is monotone, the resulting scheme is consistent, monotone and $L^\infty$ stable, and therefore convergent via the monotone convergence theorem \cite{S85}. Typical  monotone examples are $\PP_1$ (piecewise linear on triangles/tetrahedra) and $\Q_1$ (piecewise multilinear on rectangles) interpolations.

\section{Numerical solution of various types of Orienteering problem} \label{s:test}

\begin{figure}[t]
\begin{center}
\begin{tabular}{c}
\hspace{-0.2cm}\includegraphics[width=0.5\textwidth]{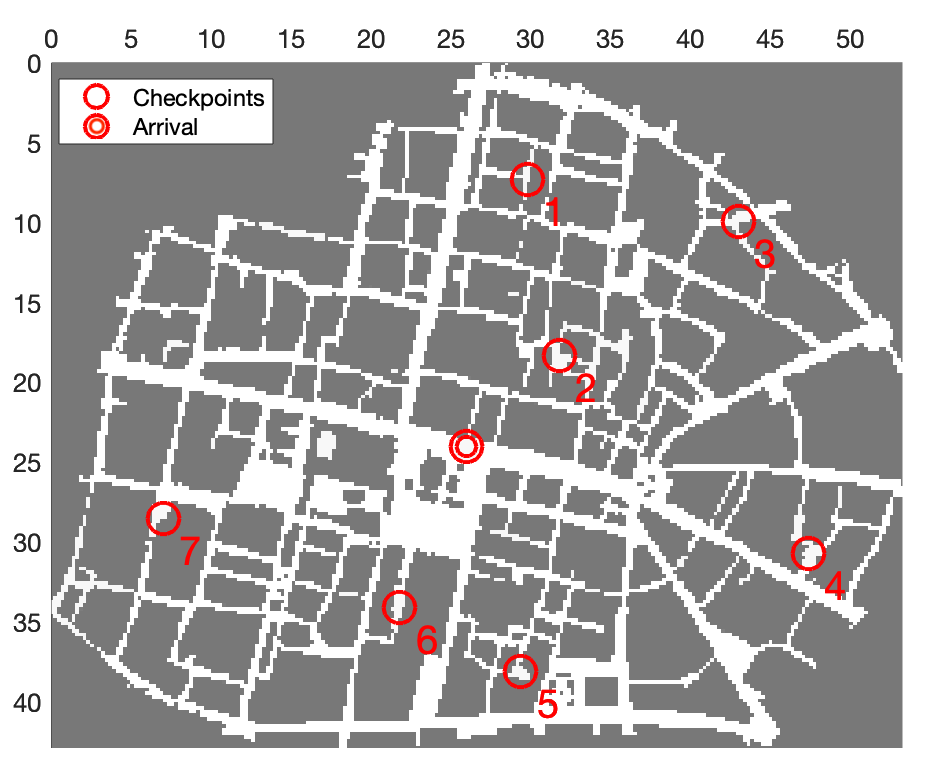}
\end{tabular}
\end{center}
\caption{Computational domain $\Omega$, ordered checkpoints and arrival.}\label{2}
\end{figure}

We consider the competition field displayed in Figure \ref{1}, i.e. the city center of Bologna, Italy. The characteristics of this urban area, basically almost flat, allow us to consider the roads as homogeneous, in the sense that all the roads, squares and alleys, can be traveled at a fixed speed, making the problem a pure minimal time problem. In order to include the obstacles and barriers we penalize the possibility to walk through them, using the running cost $\ell$. Practically speaking, we deduce from the map of Fig. \ref{1} a domain $\Omega\subset \R^2$ and an obstacle set $\Lambda\subset \Omega$. Thus we define the running cost function as 
\begin{equation}\label{runn}
 \ell(x)=\left\{ 
 \begin{array}{ll}
    1 & x\in \Omega\setminus \Lambda\\
    1/\varepsilon & x\in \Lambda.
 \end{array}\right.
\end{equation}
The dynamics of the player are isotropic, and of constant speed, while the discount factor is very small, i.e.
$$  f(x,\alpha,q)=\alpha, \quad \lambda=10^{-5}. $$
The control set is chosen as $A=B(0,1):=\{x\in\R^2 \hbox{ s.t. }\|x\|\leq 1\}$.
We set $m=7$ checkpoints $\{\cT_j\}$ on $\Omega$ as well as the arrival $\cA$ as reported in Fig. \ref{2}. Due to the nature of the feedback control that technique proposed is able to obtain, we do not need to include any information about the start since the problem is solved for any possible starting point. 

\subsection{The cross-country competition}
As stated before, in the cross-county type of competition, the checkmarks must be visited in a fixed order. This fact makes the set of discrete states $\cI$ considerably smaller, since, if for simplicity we consider that the order of the checkmarks  is from $1$ to $m$, the set $\cI$ is formed by $m+1$ elements where the $j$-element is 
$$ q_j:=(\underbrace{1,\ldots,1}_\text{$m-j-1$},\underbrace{0,\ldots,0}_\text{$j-1$}).$$ We fix the penalization parameter $\varepsilon=10^-3$.

We discretize the domain, formally considered in $[0,43]\times[0,54]$ (scale $1:12.5$) with constant space step $\Delta x=0.1$, obtaining a computational domain of $429\times 533\times 8\approx 2\cdot 10^6$ nodes. We implement the scheme \eqref{Scheme1} in semiLagrangian form. The numerical solution is obtained via a fixed point iteration, i.e.
\begin{equation}\label{sc:t}
V^{n+1}_\Delta(x_j,q) = \min \left( NV^n_\Delta(x_j,q), \Sigma\left(x_j,q,V^n_\Delta\right)\right).
\end{equation}
The equilibrium is reached within a tolerance of $10^{-3} $ in about $10^4$ iterations. 
 
\begin{figure}[t]
\begin{center}
\begin{tabular}{c}
\hspace{-0.2cm}\includegraphics[width=0.5\textwidth]{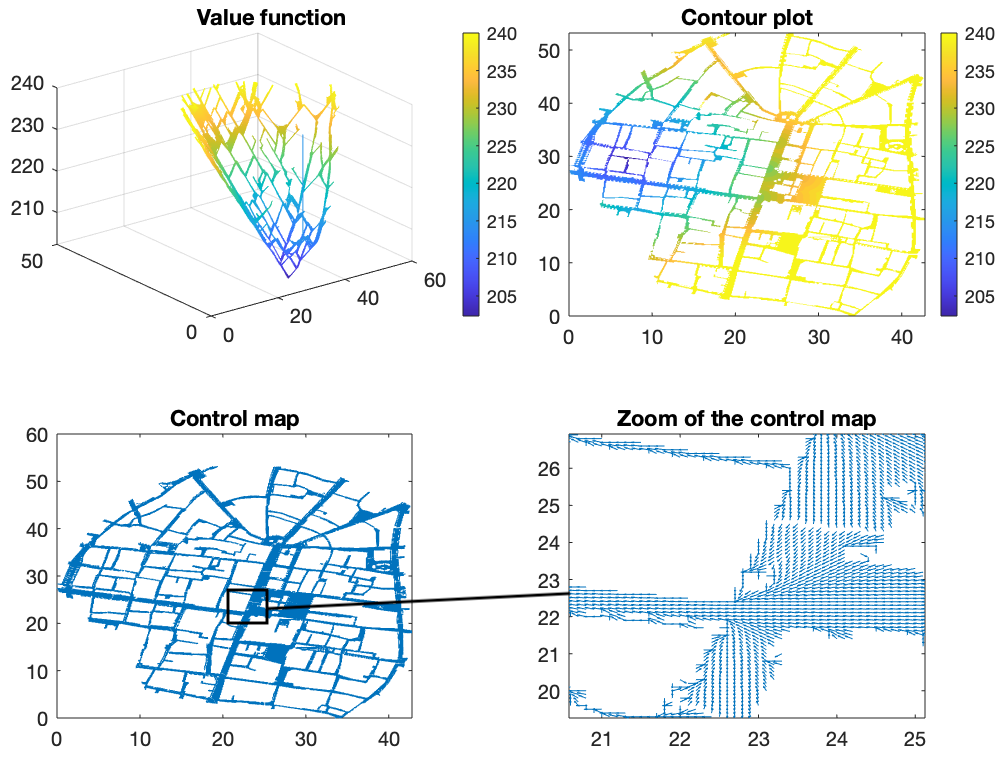}
\end{tabular}
\end{center}
\caption{Value function and control map of the state $8\cong(0,0,0,0,0,0,0)$}\label{3}
\end{figure}

 In Figure \ref{3} are shown the value function and the control mapping in the state $\hat q=8 \cong (0,0,0,0,0,0,0)$. In this case, because of the order of visit of the checkpoints, the control map leads to the only checkpoint one. Differently from a simpler minimal time problem for reaching such a point, the minimum of the value function is bigger than zero, since it contains some information about all the other points still to visit. In the zoom of the control map of the same figure is possible to observe an area where the control diverges between two parallel roads leading to the target. 

\begin{figure}[t]
\begin{center}
\begin{tabular}{c}
\hspace{-0.2cm}\includegraphics[width=0.5\textwidth]{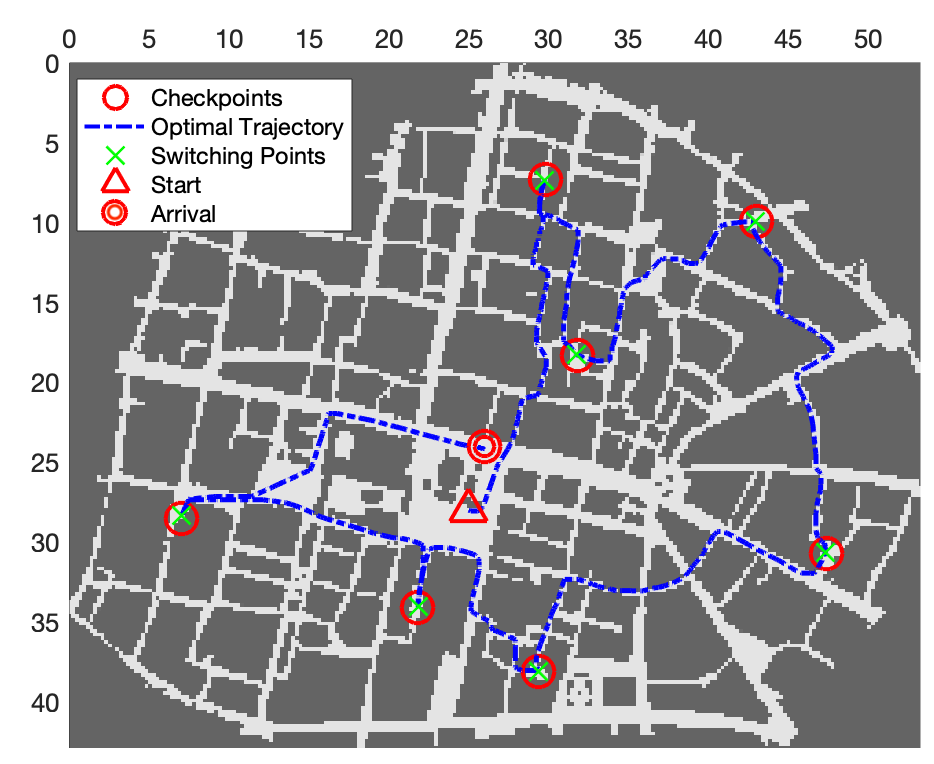}
\end{tabular}
\end{center}
\caption{Solution of the cross-country-type race. The order of the checkpoints is fixed, all the checkpoints must be visited before finishing.}\label{4}
\end{figure}

 In Figure \ref{4} we can observe the optimal trajectory of the whole circuit: the player starts in $\hat x=(28,25)$, visits all the checkpoints and then finishes reaching the arrival. The figure highlights also the switching points that are, as consequence of the choice of an adequate penalization parameter $\varepsilon$ on the checkpoints. 
Scaling the speed of the competitor to an average of $8$ Km/h,  it needs ($\tau(\hat x,\hat q)\approx$) $15.48$ minutes to conclude the course, running a total of $2.06$ Km.

\subsection{The free-order competition}
We keep the position of start, checkpoints considering a free-order race. The only difference of the free order of visit of the checkpoints requires a much larger number of states. The set $\cI$ is composed by all the possible permutations with repetitions of the elements $\{0,1\}$ on a  $7$-elements vector. Therefore, the total number of the discrete states rises to $2^m=128$, bringing the number of nodes in the discrete space to $3\cdot 10^7$. This fact, together with the number of fixed point iterations necessary to reconstruct the solution, requires a certain care from the point of view of the computational efficiency. In any case, our code, implemented in Matlab2019b on a 1,4 GHz Quad-Core Intel Core i5 processor, reaches the solution under a prescribed tolerance of $10^{-3}$ on the uniform norm in about $20$ hours of computations. If more effective performances are needed, a viable strategy would be to adapt a fast marching technique in this case (see the classic book by \cite{sethian1999fast} for details).

\begin{figure}[t]
\begin{center}
\begin{tabular}{c}
\hspace{-0.2cm}\includegraphics[width=0.5\textwidth]{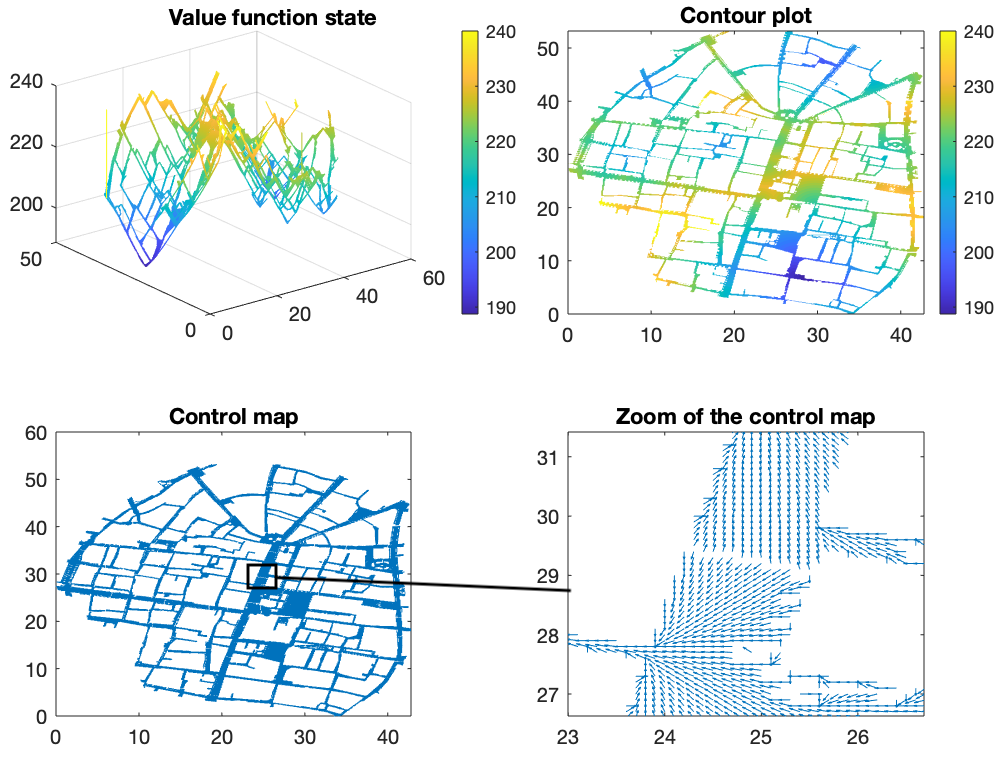}
\end{tabular}
\end{center}
\caption{Value function and control map of the state $128\cong(0,0,0,0,0,0,0)$}\label{5}
\end{figure}

Figure \ref{5} contains the value function and the control mapping of the first state $\hat q=128 \cong (0,0,0,0,0,0,0)$. It is possible to observe how, differently from the cross-county case, the minimum of the function corresponds to various points of the domain as the control map, which leads to different checkpoints depending on the current state of the player.

\begin{figure}[t]
\begin{center}
\begin{tabular}{c}
\includegraphics[width=0.5\textwidth]{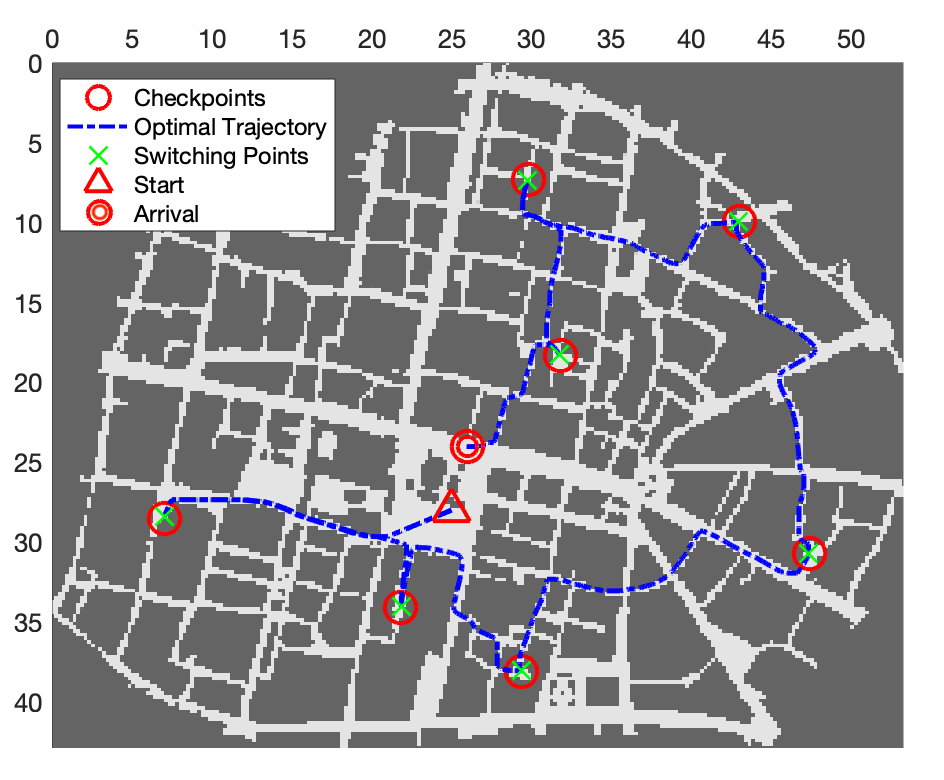}
\end{tabular}
\end{center}
\caption{Solution of the free-type race. The ordered followed, if marked as in Fig. \ref{1} is $[7,6,5,4,3,1,2]$.}\label{6}
\end{figure}

The optimal trajectory of the player is shown in Figure \ref{6}. In this case, the optimal order of visiting of the checkpoints is $[7,6,5,4,3,1,2]$. It is interesting to see how, in this case, due to the conformation of the domain, the optimal path includes various roads going through and back from the same way. The competitor, again moving to an average speed of $8$ Km/h, needs $13.92$ minutes to conclude the course, running a total of $1.85$ Km. Clearly, not bounded to a fixed order of visiting of the points, the length and consequently the time necessary to conclude the race is lower.

\subsection{The score competition}
In the score race, a competitor has at his disposition a fixed time $T>0$ to reach as many checkpoints as possible and to finish at the arrival $\cA$. To simplify this case, we assume the same value for any checkpoint. A generalization of this case is always possible, once again, through an opportune modification of the switching cost function $C$. We solve the time-dependent HJB equation \eqref{hjbt} using the same scheme \eqref{sc:t} with the only difference that now the iteration has the meaning of a time iteration and $h$ in \eqref{Scheme2} must be substituted by $\Delta t$. A clear advantage of the semiLagrangian scheme is that any CFL condition is needed to ensure stability of the scheme. This allows us to choose a large $\Delta t=5\Delta x$ and reduce in this way the dimension of the unknown $V_\Delta$. This is very important since the memory occupation is the real bottleneck in this case. 

\begin{figure}[t]
\begin{center}
\begin{tabular}{c}
\includegraphics[width=0.5\textwidth]{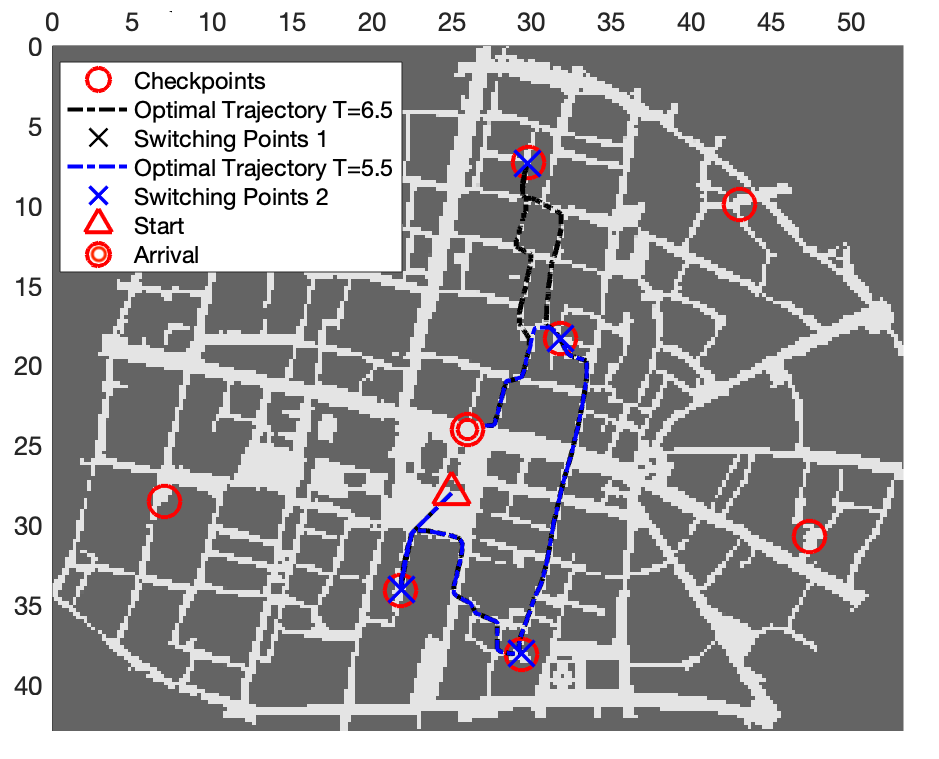}
\end{tabular}
\end{center}
\caption{Solution of the score-type race. The time at disposition for the race are 5.5 (blue line) and 6.5 (black line). }\label{7}
\end{figure}

We solve the problem in a time horizon $T=7$ minutes. Once having done it, we can compute the optimal trajectories of any case where the player has $T'<T$ as race time, just changing the initial point of the trajectory. In Figure \ref{7} we show two optimal solutions, with $T'=5.5$, $6.5$. We can see as in the first case the player choose to reach only the points $[6, 5, 2]$ before returning to the arrival, while in the second case, with a short deviation, permitted by the larger time at its disposition, it's able to reach also the point $1$ before finishing.

\section{Conclusion}\label{s:con}
We have proposed a hybrid control-based technique for the orienteering problem. Both the navigation through an inhomogeneous media (race environment/terrain) and the order/number of checkpoints' strategic choice to reach is taken into account. Many possible extensions are also possible. The most interesting is to introduce a more complex running cost $\ell$. For example, in a mountain or rural area outside the city, a critical point in the navigation choices is the altitude of the ground. The slopes and the dips reported in the map by some contour lines affect the optimal trajectories' options considerably, which may be substantially longer to avoid steeper climbs or damaging up and downs. 

Another point that deserves a more proper investigation is related to computational complexity. At this stage, the offline work of analysis of the problem is considerably high and limits the technique's real applicability. An improvement in this direction would be extremely beneficial. 

On the theoretical side, a deeper understanding of the case where $\epsilon\rightarrow 0^+$ could be interesting. We guess that the solution of the relaxed problem should converge to a solution of a problem with discontinuous running and switching cost, but such a problem may show non uniqueness of the solution. Therefore an appropriate characterization of the right solution to choose could be necessary.

\end{document}